\documentclass{article}
\begin{document}
\font\germ=eufm10
\def\ssl{\hbox{\germ sl}}
\def\slh{\widehat{\ssl_2}}
\def\ge{\hbox{\germ g}}

\makeatletter
\def\aaa{@}
\centerline{}
\centerline{\Large\bf Extremal Projectors of $q$-Boson Algebras}
\vskip7pt
\vskip15pt
\centerline{NAKASHIMA Toshiki}
\vskip10pt
\centerline{Department of Mathematics,}
\centerline{Sophia University, Tokyo 102-8554, JAPAN}
\centerline{e-mail:\,\,toshiki@mm.sophia.ac.jp}
\vskip10pt
\centerline{Abstract}
We define the extremal projector
of the $q$-boson Kashiwara algebra $B_q(\ge)$ and study
their basic properties. Applying their proerties to
 the representation theory of the category ${\cal O}(B_q(\ge))$,
whose objects are ''upper bounded '' $B_q(\ge)$-modules,
we obtain its semi-simplicity and the classification of simple modules.

\makeatother

\renewcommand{\labelenumi}{$($\roman{enumi}$)$}
\renewcommand{\labelenumii}{$(${\rm \alph{enumii}}$)$}
\font\germ=eufm10

\def\AA{{\cal A}}
\def\al{\alpha}
\def\bq{B_q(\ge)}
\def\bqm{B_q^-(\ge)}
\def\bqz{B_q^0(\ge)}
\def\bqp{B_q^+(\ge)}
\def\beneme{\begin{enumerate}}
\def\beq{\begin{equation}}
\def\beqn{\begin{eqnarray}}
\def\beqnn{\begin{eqnarray*}}
\def\bigsl{{\hbox{\fontD \char'54}}}
\def\bbra#1,#2,#3{\left\{\begin{array}{c}\hspace{-5pt}
#1;#2\\ \hspace{-5pt}#3\end{array}\hspace{-5pt}\right\}}
\def\cd{\cdots}
\def\CC{\hbox{\bf C}}
\def\ddd{\hbox{\germ D}}
\def\del{\delta}
\def\Del{\Delta}
\def\Delr{\Delta^{(r)}}
\def\Dell{\Delta^{(l)}}
\def\Delb{\Delta^{(b)}}
\def\Deli{\Delta^{(i)}}
\def\ei{e_i}
\def\eit{\tilde{e}_i}
\def\eneme{\end{enumerate}}
\def\ep{\epsilon}
\def\eeq{\end{equation}}
\def\eeqn{\end{eqnarray}}
\def\eeqnn{\end{eqnarray*}}
\def\fit{\tilde{f}_i}
\def\FF{{\rm F}}
\def\ft{\tilde{f}}
\def\gau#1,#2{\left[\begin{array}{c}\hspace{-5pt}#1\\
\hspace{-5pt}#2\end{array}\hspace{-5pt}\right]}
\def\ge{\hbox{\germ g}}
\def\gl{\hbox{\germ gl}}
\def\hom{{\hbox{Hom}}}
\def\ify{\infty}
\def\io{\iota}
\def\kp{k^{(+)}}
\def\km{k^{(-)}}
\def\llra{\relbar\joinrel\relbar\joinrel\relbar\joinrel\rightarrow}
\def\lan{\langle}
\def\lar{\longrightarrow}
\def\lm{\lambda}
\def\Lm{\Lambda}
\def\mapright#1{\smash{\mathop{\longrightarrow}\limits^{#1}}}
\def\mm{{\bf{\rm m}}}
\def\nd{\noindent}
\def\nn{\nonumber}
\def\nnn{\hbox{\germ n}}
\def\catob{{\cal O}(B)}
\def\oint{{\cal O}_{\rm int}(\ge)}
\def\ot{\otimes}
\def\op{\oplus}
\def\opi{\ovl\pi_{\lm}}
\def\ovl{\overline}
\def\plm{\Psi^{(\lm)}_{\io}}
\def\qq{\qquad}
\def\q{\quad}
\def\qed{\hfill\framebox[2mm]{}}
\def\QQ{\hbox{\bf Q}}
\def\qi{q_i}
\def\qii{q_i^{-1}}
\def\ran{\rangle}
\def\rlm{r_{\lm}}
\def\ssl{\hbox{\germ sl}}
\def\slh{\widehat{\ssl_2}}
\def\ti{t_i}
\def\tii{t_i^{-1}}
\def\til{\tilde}
\def\tt{{\hbox{\germ{t}}}}
\def\ttt{\hbox{\germ t}}
\def\ua{U_{\AA}}
\def\ue{U_{\vep}}
\def\uq{U_q(\ge)}
\def\ufin{U^{\rm fin}_{\vep}}
\def\ufinp{(U^{\rm fin}_{\vep})^+}
\def\ufinm{(U^{\rm fin}_{\vep})^-}
\def\ufinz{(U^{\rm fin}_{\vep})^0}
\def\uqm{U^-_q(\ge)}
\def\uqp{U^+_q(\ge)}
\def\uqmq{{U^-_q(\ge)}_{\bf Q}}
\def\uqpm{U^{\pm}_q(\ge)}
\def\uqq{U_{\bf Q}^-(\ge)}
\def\uqz{U^-_{\bf Z}(\ge)}
\def\ures{U^{\rm res}_{\AA}}
\def\urese{U^{\rm res}_{\vep}}
\def\uresez{U^{\rm res}_{\vep,\ZZ}}
\def\util{\widetilde\uq}
\def\uup{U^{\geq}}
\def\ulow{U^{\leq}}
\def\bup{B^{\geq}}
\def\blow{\ovl B^{\leq}}
\def\vep{\varepsilon}
\def\vp{\varphi}
\def\vpi{\varphi^{-1}}
\def\VV{{\cal V}}
\def\xii{\xi^{(i)}}
\def\Xiioi{\Xi_{\io}^{(i)}}
\def\WW{{\cal W}}
\def\wtil{\widetilde}
\def\what{\widehat}
\def\wpi{\widehat\pi_{\lm}}
\def\ZZ{\hbox{\bf Z}}

\renewcommand{\thesection}{\arabic{section}}
\section{Introduction}
\setcounter{equation}{0}
\renewcommand{\theequation}{\thesection.\arabic{equation}}

In \cite{N}, we studied the so-called $q$-boson Kashiwara algebra,
in particular, a kind of the $q$-vertex operators and their
2 point functions. We found therein some interesting object
$\Gamma$. But at that time we did not reveal its whole properties, as
``Extremal Projectors''.
Tolstoy,V.N., {\it et.al.,} introduced the notion of
``Extremal Projectors''for Lie (super)algebras and quantum (super)
algebras, and made extensive study of their properties and applied
it to the representation theory,
(see \cite{KT},\cite{To} and the references therein).
In the present paper, we shall re-define the
extremal projector for the $q$-boson algebras,
clarify their properties and apply it to the representation theory
of $q$-boson algebras.

To be more precise, let $\{e''_i, f_i, q^h\,|\,i\in I, \,\,h\in P^*\}$
be the generators of the $q$-boson algebra $\bq$.
The extremal projector $\Gamma$ is an element in
$\what\bq$(some completion of $\bq$) which satisfies the following;
\beqnn
&e''_i\Gamma=\Gamma f_i=0, \q \Gamma^2=\Gamma,\\
& \sum_k a_k\Gamma b_k=1,
\eeqnn
for some $a_k\in \bqp$ and $b_k\in \bqm$
(see Theorem \ref{ext-proj}).
Let ${\cal O}(B)$ be the category of `upper bounded' $\bq$-modules
(see Sect 3).
By using the above properties of $\Gamma$, we shall show that
the category ${\cal O}(B)$  is semi-simple and
classify its simple modules.

In \cite{K}, Kashiwara gave the projector $P$
for $q$-boson algbera of $\ssl_2$-case in order to define
the crystal base of $\uqm$. He uses it to show the
semi-simplicity of ${\cal O}(B_q(\ssl_2))$.
So our $\Gamma$ is a generalizations of his
projector $P$ to arbitrary Kac-Moody algebras.

The organization of this article is as follows;
In Sect.2, we review the definitions of
the quantum algebras and the $q$-boson Kashiwara
algebras and their properties.
In Sect.3, we introduce the category of modules of the $q$-boson
algebras ${\cal O}(B)$, which we treat in the sequel.
In Sect.4, we review so-called Drinfeld Killing form and by using it
we define some element ${\cal C}$
in the tensor product of $q$-boson algebras, which plays
a significant role of studying extremal projectors.
In Sect.5, we define extremal projectors for the $q$-boson algebras and
involve their important properties.
In the last section, we apply it to show the semi-simplicity
 of the category ${\cal O}(B)$ and
classify the simple modules in ${\cal O}(B)$.
In \cite{N}
we gave the proof of its semisimplicity, but there was a quite big gap.
Thus, the last section would be devoted to an erratum for it.
We can find an elementary proof of the semi-simplicity of the
category ${\cal O}(B)$ in {\it e.g.}\cite{Tan}.

The author would like to acknowledge Y.Koga for valuable discussions and
A.N.Kirillov for introducing the papers \cite{KT},\cite{T} to him.

\renewcommand{\thesection}{\arabic{section}}
\section{Quantum algebras and $q$-boson Kashiwara algebras}
\setcounter{equation}{0}
\renewcommand{\theequation}{\thesection.\arabic{equation}}

\newtheorem{pro2}{Proposition}[section]
\newtheorem{df2}[pro2]{Definition}
\newtheorem{lem2}[pro2]{Lemma}
\newtheorem{thm2}[pro2]{Theorem}

We shall define the algebras playing a significant role in this paper.
First, let $\ge$ be a  symmetrizable Kac-Moody algebra over {\bf Q}
with a Cartan subalgebra $\ttt$, $\{\al_i\in\ttt^*\}_{i\in I}$
 the set of simple roots and
$\{h_i\in\ttt\}_{i\in I}$  the set of coroots,
where $I$ is a finite index set. We define an inner product on
$\ttt^*$ such that $(\al_i,\al_i)\in{\bf Z}_{\geq 0}$ and
$\lan h_i,\lm\ran=2(\al_i,\lm)/(\al_i,\al_i)$ for $\lm\in\ttt^*$.
Set $Q=\oplus_i\ZZ\al_i$, $Q_+=\oplus_i\ZZ_{\geq0}\al_i$ and
$Q_-=-Q_+$. We call $Q$ a root lattice.
Let  $P$  a lattice of $\ttt^*$ {\it i.e.} a free
{\bf Z}-submodule of $\ttt^*$ such that
$\ttt^*\cong {\hbox{\bf Q}}\ot_{\ZZ}P$,
and $P^*=\{h\in \ttt|\lan h,P\ran\subset\ZZ\}$.
Now, we introduce the symbols
$\{e_i,e''_i,f_i,f'_i\,(i\in I),q^h\,(h\in P^*)\}$.
These symbols satisfy the following relations;
\beqn
&&q^0=1, \q{\hbox{\rm and }}\q q^hq^{h'}=q^{h+h'},\label{2.1}\\
&& q^he_iq^{-h}=q^{\lan h,\al_i\ran}e_i,\label{2.2}\\
&&q^he''_iq^{-h}=q^{\lan h,\al_i\ran}e''_i,\label{2.3}\\
&&q^hf_iq^{-h}=q^{-\lan h,\al_i\ran}f_i,\label{2.4}\\
&&q^hf'_iq^{-h}=q^{-\lan h,\al_i\ran}f'_i,\label{2.5}\\
&&[e_i,f_j]=\del_{i,j}(t_i-t^{-1}_i)/(q_i-q^{-1}_i),\label{2.6}\\
&&e''_if_j=q_i^{{\lan h_i,\al_j\ran}}f_je''_i+\del_{i,j},\label{2.7}\\
&&f'_ie_j=q_i^{{\lan h_i,\al_j\ran}}e_jf'_i+\del_{i,j},\label{2.8}\\
&& \sum_{k=0}^{1-{\lan h_i,\al_j\ran}}
(-1)^kX_i^{(k)}X_jX_i^{(1-{\lan h_i,\al_j\ran}-k)}=0,
\,(i\ne j),\qq\qq\qq \label{2.9}\\
&&{\hbox{ for }} X_i=e_i,\,e''_i,\,f_i,\,f'_i.\nn
\eeqn
where $q$ is transcendental over $\QQ$ and
we set $q_i=q^{(\al_i,\al_i)/2}$, $t_i=q_i^{h_i}$,
$[n]_i=(q^n_i-q^{-n}_i)/(q_i-q_i^{-1})$, $[n]_i!=\prod_{k=1}^n[k]_i$
and $X_i^{(n)}=X_i^n/[n]_i!$.

Now, we define the
algebras $B_q(\ge)$, $\ovl B_q(\ge)$ and $U_q(\ge)$.
The algebra $B_q(\ge)$ (resp. $\ovl B_q(\ge)$)
is an associative algebra generated by the symbols $\{e''_i,f_i\}_{i\in I}$
(resp. $\{e_i,f'_i\}_{i\in I}$) and $q^h$ ($h\in P^*$) with
the defining relations (\ref{2.1}), (\ref{2.3}), (\ref{2.4}), (\ref{2.7})
and (\ref{2.9})
(resp. (\ref{2.1}), (\ref{2.2}), (\ref{2.5}), (\ref{2.8}) and (\ref{2.9})) over $\QQ(q)$.
The algebra $U_q(\ge)$ is the usual quantum algebra
generated by the symbols $\{e_i,f_i\}_{i\in I}$ and $q^h$ ($h\in P^*$)
with the defining relations (\ref{2.1}),(\ref{2.2}),(\ref{2.4}),
(\ref{2.6}) and (\ref{2.9}) over $\QQ(q)$.
We shall call algebras $B_q(\ge)$ and $\ovl B_q(\ge)$ the
$q$-{\it boson Kashiwara algebras}
([K1]).
Furthermore, we define their subalgebras
\beqnn
&&T=\lan q^h|h\in P^*\ran=B_q(\ge)\cap\ovl B_q(\ge)\cap U_q(\ge),\\
&&B^{\vee}_q(\ge)\,({\hbox{resp. }}\ovl B^{\vee}_q(\ge))=
\lan e''_i,\,f_i\,({\hbox{resp. }}e_i,\,f'_i)|i\in I\ran
\subset B_q(\ge)\,({\hbox{resp. }}\ovl B_q(\ge)),\\
&&U^+_q(\ge)\,({\hbox{resp. }}U^-_q(\ge))=
\lan e_i\,({\hbox{resp. }}f_i)|i\in I\ran=:\ovl B^+_q(\ge)
\,({\hbox{resp. }}\bqm),\\
&&\uup_q(\ge)\,({\hbox{resp. }}\ulow_q(\ge))=
\lan e_i\,({\hbox{resp. }}f_i),q^h|i\in I,\,h\in P^*\ran,\\
&&B^+_q(\ge)\,({\hbox{resp. }}\ovl B^{\,\,-}_q(\ge))=
\lan e''_i\,({\hbox{resp. }}f'_i)|i\in I\ran\subset B^{\vee}_q(\ge)
\,({\hbox{resp. }}\ovl B^{\vee}_q(\ge)),\\
&&\bup_q(\ge)\,({\hbox{resp. }}\blow_q(\ge))=
\lan e''_i\,({\hbox{resp. }}f'_i),q^h|i\in I,\,h\in P^*\ran
\subset B_q(\ge)\,({\hbox{resp. }}\ovl B_q(\ge)).\\
\eeqnn
We shall use the abbreviated notations
$U$, $B$, $\ovl B$, $B^{\vee}$,$\cdots$
for $\uq$, $B_q(\ge)$, $\ovl B_q(\ge)$, $B^{\vee}_q(\ge)$,$\cdots$
if there is no confusion.

For $\beta=\sum m_i\al_i\in Q_+$ we set $|\beta|=\sum m_i$ and
$$
U^{\pm}_{\pm\beta}=\{u\in U^{\pm}|q^huq^{-h}
=q^{\pm\lan h,\beta\ran}u\,\,(h\in P^*)\},
$$
and call $|\beta|$ a height of $\beta$ and
$U^+_{\beta}$ (resp. $U^-_{-\beta}$) a weight space of
$U^+$ (resp. $U^-$) with a weight $\beta$ (resp. $-\beta$). We also define
$B^+_{\beta}$ and $\ovl B^{\,\,-}_{-\beta}$ by the similar manner.

\begin{pro2}[\cite{N}]
\begin{enumerate}
\item
We have the following algebra homomorphisms :
$\Del:U\lar U\ot U$,
$\Delr: B\lar B\ot U$, $\Dell:\ovl B\lar U\ot\ovl B$
and $\Delb:U\lar\ovl B\ot B$ given by
\beqn
&&\hspace{-30pt} \Del(q^h)=\Delr(q^h)=\Dell(q^h)=\Delb(q^h)=q^h\ot q^h,\\
&&\hspace{-30pt} \Del(e_i)=e_i\ot 1+\ti\ot e_i,\q
\Del(f_i)=f_i\ot\tii+1\ot f_i,\\
&&\hspace{-30pt} \Delr(e''_i)=(\qi-\qii)\ot\tii e_i+e''_i\ot\tii,\q
\Delr(f_i)=f_i\ot\tii+1\ot f_i,\qq\\
&&\hspace{-30pt} \Dell(e_i)=e_i\ot 1+\ti\ot e_i,\q
\Dell(f'_i)=(\qi-\qii)\ti f_i\ot 1+\ti\ot f'_i,\qq\\
&&\hspace{-30pt} \Delb(e_i)=\ti\ot{{\ti e''_i}\over{\qi-\qii}}+e_i\ot 1,\q
\Delb(f_i)=1\ot f_i+{{\tii f'_i}\over{\qi-\qii}}\ot\tii,\qq
\eeqn
and extending these to the whole algebras
by the rule: $\Del(xy)=\Del(x)\Del(y)$
and $\Deli(xy)=\Deli(x)\Deli(y)$ ($i=r,l,b$).
\item
We have the following anti-isomorphisms $S:U\lar U$ and
$\vp:\ovl B\lar B$ given by
\beqnn
&&S(e_i)=-\tii e_i,\qq S(f_i)=-f_i\ti,\qq S(q^h)=q^{-h}\\
&&\vp(e_i)=-{1\over{\qi-\qii}}e''_i,\qq
\vp(f'_i)=-(\qi-\qii)f_i,\qq \vp(q^h)=q^{-h},\\
\eeqnn
and extending these to the whole algebras
by the rule: $S(xy)=S(y)S(x)$ and $\vp(xy)=\vp(y)\vp(x)$.
Here $S$ is called a anti-pode of $U$.
We also denote $\vp|_{\uup}=\vp|_{\bup}$ by $\vp$.
\end{enumerate}
\end{pro2}

 We obtain the following triangular decompostion of
the $q$-boson Kashiwara algebra;
\begin{pro2}
\label{tri}
The multiplication map defines an isomorphism of vector spaces:
$$
\begin{array}{cccc}
&\bqm\ot T \ot \bqp&\mapright\sim&\bq \\
&u_1\ot u_2\ot u_3 & \mapsto &u_1 u_2 u_3.
\end{array}
$$
\end{pro2}

{\sl Proof.}
By \cite[(3.1.2)]{K}, we have
$$
{e''_i}^nf_j^{(m)}=
\left\{
\begin{array}{ll}
{\displaystyle \sum_{i=0}^{{\rm min}(n,m)}}
q_i^{2nm+(n+m)i-i(i+1)/2}\gau{\scriptstyle{\rm min}(n,m)},{\scriptstyle i}
f_i^{(m-i)}{e''_i}^{n-i},&{\rm if }\,\, i=j,\\
q_i^{nm\lan h_i,\al_j\ran}f_j^{(m)}{e''_i}^n,&{\rm otherwise.}
\end{array}\right.
$$
By this formula and the standard argument, we can show the proposition.\qed

We define weight completions of $L^{(1)}\ot\cdots\ot L^{(m)}$,
where $L^{(i)}=B$ or $U$.(See[T])
$$
\widehat L^{(1)}\widehat\ot\cdots\widehat\ot\widehat L^{(m)}=
\lim_{\longleftarrow\atop l}L^{(1)}\ot\cdots\ot L^{(m)}
  /(L^{(1)}\ot\cdots\ot L^{(m)})L^{+,l},
$$
where
$L^{+,l}=\oplus_{|\beta_1|+\cdots+|\beta_m|\geq l}
{L^{(1)}}^+_{\beta_1}\ot\cdots\ot {L^{(m)}}^+_{\beta_m}$.  (Note that
 $U\cong U^-\ot T\ot U^+$ and $B\cong B^-\ot T\ot B^+$. )
The linear maps $\Del$, $\Delr$, $S$, $\vp$, multiplication,
{\it e.t.c.} are
naturally extend for such completions.

\renewcommand{\thesection}{\arabic{section}}
\section{Category $\catob$}
\setcounter{equation}{0}
\renewcommand{\theequation}{\thesection.\arabic{equation}}

Let $\catob$  be the category of left
$B$-modules such that
\begin{enumerate}
\item
Any object $M$ has a weight space decomposition
$M=\oplus_{\lm\in P}M_{\lm}$ where
$M_{\lm}=\{u\in M\,|\,q^h u=q^{\lan h,\lm\ran}{\hbox{ for any }}h\in P^*\}$.
\item
For any element $u\in M$ there exists $l>0$ such that
$e''_{i_1}e''_{i_2}\cdots e''_{i_l}u=0$
 for any $i_1,i_2,\cdots,i_l\in I$.
\end{enumerate}

The similar category ${\cal O}(B^{\vee})$
for $\bq^{\vee}$ is introduced in \cite{K}, which is defined
with the above condition (ii).
In \cite{K}, Kashiwara mentions that the category
${\cal O}(B^{\vee})$ is semi-simple though he does not give an
exact proof. Here we give a proof of the semi-simplicity of $\catob$
in Sect 6.

Here for $\lm\in P$ we define the $B$-module $H(\lm)$ by
$H(\lm):=B/I_{\lm}$ where the left ideal $I_{\lm}$ is defined as 
$$
I_{\lm}:=\sum _i Be''_i+\sum_{h\in P^*} B(q^h-q^{\lan h,\lm\ran}).
$$
In Sect.6, we shall also show that $\{H(\lm)|\lm\in P\}$ is the
isomorphism class of simple modules.

\renewcommand{\thesection}{\arabic{section}}
\section{Bilinear forms and elements ${\cal C}$}
\setcounter{equation}{0}
\renewcommand{\theequation}{\thesection.\arabic{equation}}
\newtheorem{pro4}{Proposition}[section]
\newtheorem{thm4}[pro4]{Theorem}
\newtheorem{lem4}[pro4]{Lemma}
\newtheorem{ex4}[pro4]{Example}

\begin{pro4}[\cite{R},\cite{T}]
\begin{enumerate}
\item
There exists the unique bilinear form
$$
\lan \q,\q\ran:\uup\times \ulow\longrightarrow \QQ(q),
$$
satisfying the following;
\beqnn
&\lan x,y_1y_2\ran=\lan \Del(x),y_1\ot y_2\ran,\q(x\in \uup,\,y_1,y_2\in \ulow),\\
&\lan x_1x_2,y\ran=\lan x_2\ot x_1,\Del(y)\ran,\q(x_1,x_2\in \uup,\,y\in \ulow),\\
&\lan q^h,q^{h'}\ran=q^{-(h|h')}\qq(h,h'\in P^*),\\
&\lan T,f_i\ran=\lan e_i,T\ran=0,\\
&\lan e_i,f_j\ran=\del_{ij}/(q_i^{-1}-q_i),
\eeqnn
where $(\,|\,)$ is an invariant bilinear form on $\ttt$.
\item
The bilinear form $\lan \,,\,\ran$ enjoys the following properties;
\beqn
&&\hspace{-30pt}\lan xq^h,yq^{h'}\ran=q^{-(h|h')}\lan x,y\ran,
\q{\hbox{for }}x\in \uup,\,y\in \ulow,
\,h,h'\in P^*,\\
&&\hspace{-30pt}{\hbox{For any }}\beta\in Q_+,\,\lan \,\,,\,\,
\ran_{|U^+_{\beta}\times U^-_{-\beta}}
{\hbox{ is non-degenerate and }}
\lan U^+_{\gamma},U^-_{-\del}\ran=0,\,{\hbox{ if }}
\gamma\ne \del.\qq\q
\eeqn
\end{enumerate}
We call this bilinear form the {\it Drinfeld-Killing form} of $U$.
\end{pro4}

For $\beta=\sum_i m_i\al_i\in Q_+ \,\,(m_i\geq 0)$, set
$k_{\beta}:=\prod_i t_i^{m_i}$, and let $\{x^{\beta}_r\}_r$  be a basis
of $U^+_{\beta}$ and $\{y^{-\beta}_r\}_r$ be the dual
basis of $U^-_{-\beta}$ with respect to the Drinfeld-Killing form.
We denote the canonical element in $U^+_{\beta}\ot U^-_{-\beta}$
with respect to the Drinfeld-Killing form
by
$$
C_{\beta}:=\sum_r x^{\beta}_r\ot y^{-\beta}_r.
$$
We set
\begin{equation}
{\cal C}:=\sum_{\beta\in Q_+}(1\ot k_{\beta}^{-1})(1\ot S^{-1})(C_{\beta})
\in U^+\what\ot U^-=U^+\what\ot B^-.
\label{C}
\end{equation}
The element ${\cal C}$ satisfies the following relations:

\begin{pro4}
\label{prop4}
\begin{enumerate}
\item
For any $i\in I$, we have
\beqn
&&
(t_i^{-1}\ot e''_i){\cal C}={\cal C}
(t_i^{-1}\ot e''_i+(q_i-q_i^{-1})t_i^{-1}e_i\ot 1),
\label{ecec}\\
&&(f_i\ot t_i^{-1}+1\ot f_i){(\vp\ot 1({\cal C}))}=
{(\vp\ot 1({\cal C}))}(f_i\ot t_i^{-1}).
\label{fcfc}
\eeqn
Here note that (\ref{ecec}) is the equation in $\uq\what\ot \bq$
and (\ref{fcfc}) is the equation in $\bq\what\ot \bq$.
\item
The element ${\cal C}$ is invertible and the inverse is given as
\beq
{\cal C}^{-1}=
\sum_{\beta\in Q_+}q^{-(\beta,\beta)}
(k_{\beta}\ot k_{\beta}^{-1})(S^{-1}\ot S^{-1})(C_{\beta})
\label{c-inverse}
\eeq
\end{enumerate}
\end{pro4}

{\sl Proof.}
The proof of (\ref{fcfc}) has been given in \cite[6.2]{N}.
Thus, let us show (\ref{ecec}). For that purpose, we need the following lemma;
\begin{lem4}
\label{lmc}
For $\beta\in Q_+$, let $C_{\beta}=\sum_r x^{\beta}_r\ot
y^{\beta}_r$ be the canonical element in
$U^+_{\beta}\ot U^-_{-\beta}$ as above and set
$C'_{\beta}:=(1\ot S^{-1})(C_{\beta})$.
Then for any $\beta\in Q_+$ and $i\in I$,
we have
\beq
[t_i^{-1}\ot e''_i, (1\ot k^{-1}_{\beta+\al_i})(C'_{\beta+\al_i})]
=(1\ot k^{-1}_{\beta})(C'_{\beta})(t^{-1}_ie_i\ot (q_i-q_i^{-1}))
\in \uq\ot \bq,
\label{ecc}
\eeq
where we use the identification $\bqm=\uqm$.
\end{lem4}

{\sl Proof.}
Applying $\lan \cdot,z\ran\ot 1$ on the both sides of
(\ref{ecc}) where $z\in U^-_{-\beta-\al_i}$, we obtain
\beqnn
(\lan\cdot,z\ran\ot 1)({\rm L.H.S. of (\ref{ecc})})
&=& \sum_r \lan t_i^{-1}x^{\beta+\al_i}_r,z\ran\ot
e''_ik^{-1}_{\beta+\al_i}
S^{-1}(y^{-\beta-\al_i}_r)\\
&&-\lan x^{\beta+\al_i}_rt_i^{-1}, z\ran
\ot k^{-1}_{\beta+\al_i}S^{-1}(y^{-\beta-\al_i}_r)e''_i\\
&=& q^{-(\al_i,\beta+\al_i)}e''_ik^{-1}_{\beta+\al_i}S^{-1}(z)
-k^{-1}_{\beta+\al_i}S^{-1}(z)e''_i\\
&=& k^{-1}_{\beta+\al_i}(e''_iS^{-1}(z)-S^{-1}(z)e''_i).
\eeqnn
\beqn
(\lan\cdot,z\ran\ot 1)({\rm R.H.S. of (\ref{ecc})})
&=& \sum_r \lan x^{\beta}_r t^{-1}_ie_i,z\ran
\ot (q_i-q_i^{-1})k^{-1}_{\beta}S^{-1}(y^{-\beta}_r)\qq
\label{zz}
\eeqn
For $z \in U^-_{-\beta-\al_i}$ we can define $v \in U^-_{-\beta}$
uniquely by
$$
\Delta(z)=1\ot z+f_i\ot vt_i^{-1}+\cd .
$$
By the property of the Drinfeld Killing form, we have
\beqnn
\lan x^{\beta}_r t^{-1}_ie_i,z\ran
&=& \lan e_i\ot x^{\beta}_r t^{-1}_i,\Delta(z)\ran\\
&=& \lan e_i\ot x^{\beta}_r t^{-1}_i, 1\ot z+f_i\ot vt_i^{-1}+\cd\ran \\
&=& \lan e_i,f_i\ran \lan x^{\beta}_r t^{-1}_i,vt_i^{-1}\ran\\
&=& \frac{q_i^{-2}}{q_i^{-1}-q_i}\lan x^{\beta}_r,v\ran
\eeqnn
Thus,
\beq
{\rm R.H.S.\,\,of\,\,(\ref{zz})}=-q_i^{-2}k^{-1}_{\beta}S^{-1}(v)
\label{ksv}
\eeq
Here in order to comlete the proof of Lemma \ref{lmc},
let us show;
\beq
e''_iS^{-1}(z)
-S^{-1}(z)e''_i=-q_i^{-2}t_iS^{-1}(v).
\label{zv}
\eeq
Without loss of generality,
we may assume that $z$ is in the form $z=f_{i_1}f_{i_2}\cd f_{i_k}\in
U^-_{-\beta-\al_i}$ ($\beta+\al_i=\al_{i_1}+\cd+\al_{i_k}$).
For $\beta=\sum_jm_j\al_j$, we shall show by the induction
on $m_i$ for fixed $i\in I$.

If $m_i=0$, $z$ is in the form
$z=z' f_iz''$ where $z'$ and $z''$ are monomials
 of $f_j$'s not including $f_i$.
By $S^{-1}(f_j)=-t_jf_j$ and $e''_i(t_jf_j)=(t_jf_j)e''_i$ $(i\ne j)$
we have
\beq
e''_iS^{-1}(z')=S^{-1}(z')e''_i,\qq
e''_iS^{-1}(z'')=S^{-1}(z'')e''_i.
\label{d-prime}
\eeq
Hence, we obtain
\beqnn
e''_iS^{-1}(z)
&=&S^{-1}(z'')(-e''_it_if_i)S^{-1}(z')
=S^{-1}(z'')(-t_if_ie''_i-q_i^{-2}t_i)S^{-1}(z')
\\
&=&
S^{-1}(z'')(-t_if_i)S^{-1}(z')e''_i
-q_i^{-2}S^{-1}(z'')t_iS^{-1}(z')\\
&=&S^{-1}(z'')S^{-1}(f_i)S^{-1}(z')e''_i
-q^{(\beta''-\al_i,\al_i)}t_iS^{-1}(z'z''),
\eeqnn
where $\beta''=wt(z'')$. Therefore, for $m_i=0$, we have
$$
{\rm L.H.S.\,\, of \,\,(\ref{zv})}=
-q^{(\beta''-\al_i,\al_i)}t_iS^{-1}(z'z'').
$$
In the case $m_i=0$ we can easily obtain
$v=q^{(\beta'',\al_i)}z'z''$ and then
$$
{\rm R.H.S.\,\, of \,\,(\ref{zv})}=
-q^{(\beta''-\al_i,\al_i)}t_iS^{-1}(z'z'')=
{\rm L.H.S.\,\, of \,\,(\ref{zv})}
$$
Thus, the case $m_i=0$ has been shown.

Suppose that $m_i>0$.
we divide $z=z'z''$ such that $m'_i<m_i$ and $m''_i< m_i$ where
$m'_i$( resp. $m''_i$) is the number of $f_i$ including in $z'$
(resp. $z''$). Writing
\beqnn
&&
\Delta(z')=1\ot z'+f_i\ot v't_i^{-1}+\cd,\\
&&\Delta(z'')=1\ot z''+f_i\ot v''t_i^{-1}+\cd,
\eeqnn
and calculating $\Del(z'z'')$ directly, we obtain
\beq
v=z'v''+q^{(\beta'',\al_i)}v'z''.
\label{vvv}
\eeq
By the hypothesis of the induction,
\beqnn
e''_iS^{-1}(z)&=&
e''_iS^{-1}(z'')S^{-1}(z')=(S^{-1}(z'')e''_i-q_i^{-2}t_iS^{-1}(v''))S^{-1}(z')\\&=& S^{-1}(z'')e''_iS^{-1}(z')-q_i^{-2}t_iS^{-1}(z'v'')\\
&=& S^{-1}(z'')(S^{-1}(z')e''_i-q_i^{-2}t_iS^{-1}(v'))
-q_i^{-2}t_iS^{-1}(z'v'')\\
&=& S^{-1}(z'z'')e''_i-
q_i^{-2}t_i(S^{-1}(z'v'')+q^{(\beta'',\al_i)}S^{-1}(v'z''))\\
&=&S^{-1}(z)e''_i-q_i^{-2}t_iS^{-1}(v).
\eeqnn
Note that in the last equality, we use (\ref{vvv}).
Now, we have completed to show Lemma \ref{lmc}. \qed

\vskip3mm
{\sl Proof of Proposition \ref{prop4}.}\,\,
If $\beta\in Q_+$ does not include $\al_i$, since
$e''_i$ and $S^{-1}(z)$ ($z\in U^-_{-\beta}$) commute with each other
by (\ref{d-prime}), we have
$$
(t_i^{-1}\ot e''_i)(1\ot k_{\beta}^{-1})(C'_{\beta})
=(1\ot k_{\beta}^{-1})(C'_{\beta})(t_i^{-1}\ot e''_i).
$$
Thus,
we have
\beqnn
(t_i^{-1}\ot e''_i){\cal C}&-&{\cal C}(t_i^{-1}\ot e''_i)\\
&=&\sum_{\gamma\in Q_+}(t_i^{-1}\ot e''_i)(1\ot k_{\gamma}^{-1})(C'_{\gamma})
-(1\ot k_{\gamma}^{-1})(C'_{\gamma})(t_i^{-1}\ot e''_i)\\
&=&\sum_{\beta\in Q_+}
(t_i^{-1}\ot e''_i)(1\ot k_{\beta+\al_i}^{-1})(C'_{\beta+\al_i})
-(1\ot k_{\beta+\al_i}^{-1})(C'_{\beta+\al_i})(t_i^{-1}\ot e''_i)\\
&=&\sum_{\beta\in Q_+}
[t_i^{-1}\ot e''_i, (1\ot k_{\beta+\al_i}^{-1})(C'_{\beta+\al_i})]\\
&=& \sum_{\beta\in Q_+}
(1\ot k_{\beta}^{-1})(C'_{\beta})((q_i-q_i^{-1})t_i^{-1}e_i\ot 1)
\qq ({\rm by \,\,Lemma \,\ref{lmc}})\\
&=&
{\cal C}((q_i-q_i^{-1})t_i^{-1}e_i\ot 1).
\eeqnn
Then we obtain (\ref{ecec}).

Next, let us show (ii).
Set
$\wtil{\cal C}:=\sum q^{(\beta,\beta)}(1\ot k_{\beta})(S\ot 1)(C_{\beta})$.
By \cite[Sect.4]{T}, we have
$\wtil{\cal C}^{-1}:=
\sum q^{(\beta,\beta)}(k^{-1}_{\beta}\ot
k_{\beta})(C_{\beta})$.
Here note that
\beqnn
(S^{-1}\ot S^{-1})(\wtil{\cal C})
&=& \sum q^{(\beta,\beta)}(1\ot S^{-1})
\{(1\ot k_{\beta})(C_{\beta})\}\\
&=& \sum q^{(\beta,\beta)}\{(1\ot S^{-1})
(C_{\beta})\}(1\ot k_{\beta}^{-1})\\
&=& \sum (1\ot k_{\beta}^{-1})(1\ot S^{-1})(C_{\beta})\\
&=&{\cal C}
\eeqnn
Thus, we obtain
\beqnn
{\cal C}^{-1}&=&
(S^{-1}\ot S^{-1})(\wtil{\cal C}^{-1})\\
&=& \sum q^{(\beta,\beta)}
\{(S^{-1}\ot S^{-1})(C_{\beta})\}(k_{\beta}\ot k_{\beta}^{-1})\\
&=&\sum q^{-(\beta,\beta)}
(k_{\beta}\ot k_{\beta}^{-1})(S^{-1}\ot S^{-1})(C_{\beta}),
\eeqnn
 and cpmpleted the proof of Proposition
\ref{prop4}. \qed

{\sl Remark.}
By the explicit form of ${\cal C}^{-1}$ in (\ref{c-inverse}), we find that
${\cal C}^{-1}\in \uqp\what\ot\uqm=\uqp\what\ot\bqm$.
\renewcommand{\thesection}{\arabic{section}}
\section{Extremal Projectors}
\setcounter{equation}{0}
\renewcommand{\theequation}{\thesection.\arabic{equation}}
\newtheorem{ex5}{Example}[section]
\newtheorem{pro5}[ex5]{Proposition}
\newtheorem{thm5}[ex5]{Theorem}
\newtheorem{lem5}[ex5]{Lemma}

Let ${\cal C}$ be as in Sect.4.
We define the {\it extremal projector } of $\bq$ by
\beq
\Gamma :=m\circ \sigma\circ(\vp\ot1)({\cal C})
=\sum_{\beta\in Q_+,\,r}k^{-1}_{\beta}S^{-1}(y^{-\beta}_r)\vp(x^{\beta}_r),
\label{gamma}
\eeq
where $m:a\ot b\mapsto ab$ is the multiplication and
$\sigma:a\ot b\mapsto b\ot a$
is the permutation.

Here note that $\Gamma$ is a well-defined element in $\what \bq$.

\begin{ex5}[\cite{K,N}]
In
$\ssl_2$-case, the following is the explicit form of $\Gamma$.
$$
\Gamma=\sum_{n\geq 0}q^{\frac{1}{2}n(n-1)}(-1)^n f^{(n)}{e''}^n.
$$
\end{ex5}

\begin{thm5}
\label{ext-proj}
The extremal projector $\Gamma$ enjoys the following properties:
\begin{enumerate}
\item
$e''_i\Gamma=0,\qq
\Gamma f_i=0\qq(\forall i\in I).$
\item
$\Gamma^2=\Gamma$.
\item
There exists $a_k\in B_q^{-}(\ge)(=\uqm)$, $b_k\in B_q^{+}(\ge)$ such that
$$
\sum_k a_k\Gamma b_k=1.
$$
\item
$\Gamma$ is a well-defined element in $\what B^{\vee}_q(\ge)$.
\end{enumerate}
\end{thm5}

{\sl Proof.}
It is easy to see (iv) by the explicit forms
of the anti-pode $S$, the anti-isomorphism
$\vp$ and $\Gamma$ in (\ref{gamma}).
The statement (ii) is an immediate consequence of (i).
So let us show (i) and (iii).
The formula $\Gamma f_i=0$ has been shown in \cite{N}. Thus,
we shall show $e''_i\Gamma=0$.
Here let us wrtie
${\cal C}=\sum_k c_k\ot d_k,$
where $c_k\in \uqp$ and $d_k\in \bqm$.
Thus, we have
$$
\Gamma=\sum_k d_k\vp(c_k).
$$
The equation (\ref{ecec}) can be written as follows;
\beq
\sum_k t_i^{-1}c_k\ot e''_i d_k
=\sum_k c_kt_i^{-1}\ot d_ke''_i
+(q_i-q_i^{-1})c_k t_i^{-1}e_i\ot d_k.
\label{cd}
\eeq
Applying $m\circ \sigma\circ(\vp\ot1)$ on the both sides of (\ref{cd}),
we get
$$
\sum_k e''_id_k\vp(c_k)t_i
=\sum_kd_ke''_it_i\vp(c_k)-\sum_k d_ke''_it_i\vp(c_k)=0,
$$
and then
$e''_i \Gamma t_i=0$, which implies the desired result
since $t_i$ is invertible.

Next, let us see (iii).
By the remark in the last section, we can write
$$
{\cal C}^{-1}=\sum_k b'_k\ot a_k\in \uqp\what\ot\bqm.
$$
Then,
\beq
1\ot 1=\sum_{j,k}b'_kc_j\ot a_k d_j.
\label{abcd}
\eeq
Applying $m\circ\sigma\circ(\vp\ot 1)$ on the both sides of (\ref{abcd}),
we obtain
$$
1=\sum_{j,k}a_kd_j\vp(c_j)\vp(b'_k)=\sum_k a_k\Gamma \vp(b'_k).
$$
Here setting $b_k:=\vp(b'_k)$, we get (iii).\qed

\renewcommand{\thesection}{\arabic{section}}
\section{Representation Theory of ${\cal O}(B)$}
\setcounter{equation}{0}
\renewcommand{\theequation}{\thesection.\arabic{equation}}
\newtheorem{thm6}{Theorem}[section]
\newtheorem{ex6}{thm6}[section]
\newtheorem{pro6}[thm6]{Proposition}
\newtheorem{lem6}[thm6]{Lemma}

As an application of the extremal projector $\Gamma$, we shall show the
following theorem;

\begin{thm6}
\label{main}
\begin{enumerate}
\item
The category $\catob$ is a semi-simple category.
\item
The module $H(\lm)$ is a simple object of $\catob$ and for any
simple ojbect $M$ in $\catob$ there exists some $\lm\in P$
such that $M\cong H(\lm)$. Furthermore, $H(\lm)$ is a
rank one free $\bqm$-module.
\end{enumerate}
\end{thm6}

In order to show this theorem, we need to prepare several things.

For an object $M$ in ${\cal O}(B)$, set
$$
K(M):=\{v\in M\,|\,e''_i v=0{\rm \,\,for \,\,any \,\,}i\in I\}.
$$
\begin{lem6}
\label{lem1}
For an object $M$ in ${\cal O}(B)$, we have
\beq
\Gamma\cdot M=K(M)
\label{gmkm}
\eeq
\end{lem6}

{\sl Proof.}
By Theorem \ref{ext-proj}(i), we have $e''_i\Gamma=0$ for any $i\in I$.
Thus, it is trivial to see that $\Gamma\cdot M\subset K(M)$.
Owing to the explicit form of $\Gamma$, we find that
$$
1-\Gamma\in \sum_i\what\bq e''_i.
$$
Therefore, for any $v\in K(M)$ we get $(1-\Gamma)v=0$, which implies that
$\Gamma\cdot M\supset K(M)$.\qed

\begin{lem6}
\label{lem2}
For an object $M$ in ${\cal O}(B)$, we have
\beq
M=\bqm\cdot(K(M))
\label{mbkm}
\eeq
\end{lem6}

{\sl Proof.}
By Theorem \ref{ext-proj}(iii), we have $1=\sum_k a_k\Gamma b_k$
($a_k\in \bqm$, $b_k\in \bqp$). For any $u\in M$,
$$
u=\sum_k a_k(\Gamma b_k u).
$$
By Lemma \ref{gmkm}, we have $\Gamma b_k u\in K(M)$.
Then we obtain the desired result.\qed

\begin{pro6}
\label{prop-oplus}
For an object $M$ in ${\cal O}(B)$, we have
\beq
M=K(M)\oplus (\sum_i{\rm Im}(f_i)).
\label{m-oplus}
\eeq
\end{pro6}

{\sl Proof.}
By (\ref{mbkm}),
we get
$$
M=K(M)+ (\sum_i{\rm Im}(f_i)).
$$
Thus, it is sufficient to show
\beq
K(M)\cap(\sum_i{\rm Im}(f_i))=\{0\}.
\label{00}
\eeq
Let $u$ be a vector in $K(M)\cap(\sum_i{\rm Im}(f_i))$.
Since $u\in \sum_i{\rm Im}(f_i)$, there exist $\{u_i\in M\}_{i\in I}$
such that
$u=\sum_{i\in I}f_iu_i$.
By the argument in the proof of Lemma \ref{lem1}, we have
$\Gamma u=u$ for $u\in K(M)$.
It follows from Theorem \ref{ext-proj}(i) that
$$
u=\Gamma u=\sum_{i\in I}(\Gamma f_i)u_i=0,
$$
which implies (\ref{00}).
\qed

\vskip3mm

\begin{lem6}
\label{fin}
If $u$, $v\in M$ ($M$ is an object in ${\cal O}(B)$) satisfies
$v=\Gamma u$, then there exists $P\in \bq$ such that
$v=P u$.
\end{lem6}

{\sl Proof.}
By the definition of the category of ${\cal O}(B)$, there exists $l>0$
such that $\vp(x^{\beta}_r)u=0$ for any $r$ and $\beta$ with
$|\beta|>l$.
Thus, by the explicit form of $\Gamma$ in (\ref{gamma}), we can write
$$
v=\Gamma u=
(\sum_{|\beta|\leq l,\,\,r}
k^{-1}_{\beta}S^{-1}(y^{-\beta}_r)\vp(x^{\beta}_r))u,
$$
which implies our desired result.\qed
\vskip3mm

{\sl Proof of Theorem \ref{main}.}
Let $L\subset M$ be objects in the category ${\cal O}(B)$.
We shall show that there exists a submodule $N\subset M$ such that
$M=L\oplus N$.

Since $K(M)$ (resp. $K(L)$) is invariant by the actin of any $q^h$, we have
the weight space decomposition;
$$
K(M)=\bigoplus_{\lm\in P}K(M)_{\lm}
\,\,({\rm resp. }\,\,
K(L)=\bigoplus_{\lm\in P}K(L)_{\lm}).
$$
There exist subspaces $N_{\lm}\subset K(M)_{\lm}$ such that
$K(M)_{\lm}=K(L)_{\lm}\oplus N_{\lm}$, which is a decomposition of
a vector space. Here set $N:=\oplus_{\lm}N_{\lm}$.
We have
$$
K(M)=K(L)\oplus N.
$$
Let us show
\beq
M=L\oplus \bq\cdot N.
\label{lmn}
\eeq
Since
$M=\bq\cdot(K(M))=\bq(K(L)\oplus N)$, we get
$M=L+\bq\cdot N$.
Let us show
\beq
L\cap \bq\cdot N=\{0\}.
\label{000}
\eeq
For $v\in L\cap \bq\cdot N$
we have by Theorem \ref{ext-proj} (iii), 
$$
v=\sum_k a_k(\Gamma _k v).
$$
It follows from  $v\in L$ that $\Gamma b_k v\in K(L)$, and
from  $v\in \bq\cdot N$ that $\Gamma b_k v\in \Gamma(\bq\cdot N)=N$.
These imply
$$
\Gamma b_k v\in K(L)\cap N=\{0\}.
$$
Hence we get $v=0$ and then (\ref{lmn}).

Next, let us show (ii).
As an immediate consequence of Proposition \ref{tri} we can see that
$H(\lm)$ is a rank one free $\bqm$-module.

Let $\pi_{\lm}:\bq\rightarrow H(\lm)$ be the canonical projection
and set $u_{\lm}:=\pi_{\lm}(1)$.
Here we have
$$
H(\lm)=\bqm\cdot u_{\lm}
=\QQ(q)u_{\lm}+\sum_i{\rm Im}(f_i).
$$
It follows from  this, Proposition \ref{prop-oplus} and
$\QQ(q)u_{\lm}\subset K(H(\lm))$ that
$H(\lm)=
\QQ(q)u_{\lm}\oplus \sum_i{\rm Im}(f_i)$ and then
\beq
\Gamma\cdot H(\lm)=K(H(\lm))=\QQ(q)u_{\lm}.
\label{ulm}
\eeq
In order to  show the irreducibility of $H(\lm)$,
it is sufficient to see that
for arbitrary $u(\ne 0),\,v\in H(\lm)$ there exists
$P\in \bq$ such that $v=P u$.
Set $v=Q u_{\lm}$ $(Q\in \bqm)$. By Theorem \ref{ext-proj} (iii),
we have
$$
u=\sum_k a_k(\Gamma b_k u)\ne 0.
$$
Then, for some $k$ we have
$\Gamma b_k u\ne0$, which implies that
$c\Gamma b_k u= u_{\lm}$
 for some non-zero scalar $c$.
Therefore, by Lemma \ref{fin}, there exists some $R\in \bq$ such that
$u_{\lm}=Ru$ and then we have
$$
v=Q u_{\lm}=QRu.
$$
Thus, $H(\lm)$ is a simple module in ${\cal O}(B)$.

Suppose that $L$ is a simple module in ${\cal O}(B)$.
First, let us show
\beq
{\rm dim}(K(L))=1.
\label{dim1}
\eeq
For $x,y(\ne 0)\in K(L)$, there exists $P\in \bq$ such that
$y=Px$. Since $x\in K(L)$, we can take $P\in \bqm$.
Because $y\in K(L)$ and $K(L)\cap \sum_i{\rm Im}(f_i)=\{0\}$,
we find that $P$ must be a scalar, say $c$.
Thus, we have $y=cx$, which derives (\ref{dim1}).

Let $u_0$ be a basis vector in $K(L)$.
The space $K(L)$ is invariant by the action of any $q^h$ and then,
$u_0\in L_{\lm}$ for some $\lm\in P$.
Therefore, since $H(\lm)$ is a rank one free
$\bqm$-module, the map
\beqnn
\phi_{\lm}&:& H(\lm)\longrightarrow L\\
&& Pu_{\lm}\mapsto P u_0,\q(P\in \bqm),
\eeqnn
is a well-defined non-trivial homomorphism of $\bq$-modules.
Thus, by Schur's lemma, we obtain
$H(\lm)\cong L$.\qed


\begin{thebibliography}{99}

\def\CMP{\sl Commum.Math.Phys.}

\def\IJMP{\sl Int.J.Mod.Phys.}

\def\Duke{\sl Duke Math.J.}

\bibitem{K} Kashiwara M., On crystal bases of the $q$-analogue of universal
    enveloping algebras, {\Duke}, {\bf 63}, 465--516, (1991).

\bibitem{KT} Khoroshkin S.M. and Tolstoy V.N., Exremal porjector and universal
R-matrix for quantized contragradient Lie (super) algebras, Quantum Groups and
related topics, edited by Gielerak et al, 23--32, (1992).

\bibitem{N}  Nakashima T.,
             Quantum $R$-matrix and Intertwiners for the Kashiwara
             algebras, {\CMP}, {\bf 164}, 239--258, (1994).

\bibitem{R} Rosso, M, Analogues de la forme de Killing et du
th\' eor\` eme d'Harish-Chandra pour les groupes quantiques,
{\sl Ann.scient.\' Ec,Norm. Sup.} {\bf 23} (1990)  445--467.

\bibitem{T} Tanisaki T, Killing forms, Harish-Chandra isomorphisms,
and universal $R$-matrices for quantum algebras,
{\IJMP} {A7} {Suppl. 1B (1992) 941--961}.

\bibitem{To} Tolstoy V.N., Extremal projectors for Quantized Kac-Moody
superalgebras and some of their applications, Quantum Groups({\it
Clausthal, 1989}), Lecture Notes in Physics {\bf 370}, Springer,
Berlin, (1990), 118--125.

\bibitem{Tan} Tan Y., The $q$-analogue of bosons and Hall algebras,
preprint.

\end{thebibliography}
\end{document}